\def\documentid{\getTime-stamp: "2004-03-09 12:34:41 shimo"}
\def\getTime-stamp: "#1"{#1}
\documentclass[12pt,leqno]{article}

\usepackage[dvips]{graphicx}
\usepackage{amsmath}
\usepackage{bibcite}

\usepackage{xr}
 \newlength{\marginL} \marginL 2.5cm
 \newlength{\marginR} \marginR 2.5cm
 \newlength{\marginU} \marginU 3cm
 \newlength{\marginD} \marginD 4cm
 \newlength{\shiftH}  \shiftH 0mm
 \newlength{\shiftV}  \shiftV 0mm

 \evensidemargin \oddsidemargin
\newcommand{\refeq}[1]{(\ref{eq:#1})}

\newcommand{\refsec}[1]{Section~\ref{sec:#1}}

\newcommand{\reflemma}[1]{Lemma~\ref{lemma:#1}}

\newtheorem{lemma}{Lemma}

\newcommand{\myproofbox}{\qquad\rule[-0.2em]{0.3em}{0.6em}}

\newcommand{\hskiptheorem}{\hskip 0.5em}

{\begin{trivlist} \item[\hskip \labelsep{\em Proof.}\hskiptheorem]}%
{\myproofbox\end{trivlist}}

{\begin{trivlist} \item[\hskip \labelsep{\it Proof of {#1}.}\hskiptheorem]}%
{\myproofbox\end{trivlist}}

\makeatletter
\def\modifytosectioncounter#1{
\@addtoreset{#1}{section}%
\expandafter\def\csname the#1\endcsname{\thesection.\arabic{#1}}}
\modifytosectioncounter{equation}
\makeatother

\newcommand{\cR}{{\cal R}}

\newcommand{\fract}[2]{\mbox{$#1 \over #2$}}

\newcommand{\hu}{\hat u}
\newcommand{\hv}{\hat v}

\newcommand{\dabdab}{(d^{ab})^2}
\newcommand{\gabpgabp}{(\phi^{abp})^2}
\newcommand{\gappgapp}{(\phi^{app})^2}

\begin{document}
\noindent{auttex04/tech.tex\hfill\documentid}

\begin{center}
 \bf\Large Technical details of \\
the multistep-multiscale bootstrap resampling\\[3ex]
 \large
 Hidetoshi Shimodaira\\[3ex]
Department of Mathematical and Computing Sciences\\
Tokyo Institute of Technology\\
2-12-1 Ookayama, Meguro-ku, Tokyo 152-8552, Japan\\
shimo@is.titech.ac.jp
\end{center}

\vspace{1cm}

\begin{center}
\bf March, 2004
\end{center}

\vspace{2cm}

\begin{center}
{\bf RUNNING HEAD: DETAILS OF MULTISCALE BOOTSTRAP} 
\end{center}

\vspace{1cm}

{\bf
\begin{quote}
Research Reports {B}-403, Department of Mathematical and Computing
  Sciences, Tokyo Institute of Technology, Tokyo, Japan.
\end{quote}

\begin{quote}
Supplement Material to 
``Approximately unbiased tests of regions using multistep-multiscale
  bootstrap resampling.''
\end{quote}
}

\clearpage
\section*{Summary}

The technical details of the new bootstrap method of Shimodaira (2004)
are given here in mathematical proofs as well as a supporting computer
program.  Approximately unbiased tests based on the bootstrap
probabilities are considered in Shimodaira (2004) for the exponential
family of distributions with unknown expectation parameter vector, where
the null hypothesis is represented as an arbitrary-shaped region with
smooth boundaries.  It has been described in the lemmas of Shimodaira
(2004) that the newly developed three-step multiscale bootstrap method
calculates an asymptotically third-order accurate $p$-value. All the
mathematical proofs of these lemmas are shown here. The straightforward,
though very tedious, calculations involving tensor notations are
verified in Shimodaira (2003), which is, in fact, a computer program for
{\it Mathematica}. Here we also give a brief explanation of this program.

\nocite{bib:Shim:2004:AUR}
\nocite{bib:Shim:2003:AAB}

\appendix
\section{Mathematical proofs of lemmas.}

The proofs of Lemma~1 to Lemma~7 of Shimodaira (2004) are presented
here.  The equation numbers indicate those of Shimodaira (2004) except
for those defined here.

\subsection{Proof of Lemma~\ref{lemma:uv-density}.} \label{sec:lemma1}

The expression of $f(y;\eta)$ is obtained from (\ref{eq:expdensity})
using the Taylor series of the reparameterization $ \theta^i \approx
\eta_i +\fract{1}{2} \phi^{ijk} \eta_j \eta_k +\fract{1}{6} \phi^{ijkl}
\eta_j \eta_k \eta_l $.  The expressions for $\psi(\theta)$ and $h(y)$
are shown in Shimodaira (2003);
\begin{align} \label{eq:thetaineta}
 \psi(\theta) &\approx \psi(0) + \fract{1}{2} (\eta_i)^2
+\fract{1}{3} \phi^{ijk} \eta_i \eta_j \eta_k
+\fract{1}{8}\phi^{ijkl} \eta_i\eta_j\eta_k\eta_l,\\
\label{eq:hiny}
  h(y) &\approx -\psi(0)+\fract{1}{2}p\log(2\pi) 
+\fract{1}{8}\phi^{iijj} -\fract{1}{6}
(\phi^{ijk})^2
-\fract{1}{2} \phi^{ijj}  y_i\\
 &+ \fract{1}{2} \left( 
\delta_{ij} + \fract{1}{2} \phi^{ikl}\phi^{jkl} - \fract{1}{2} \phi^{ijkk}
\right) y_i y_j 
+\fract{1}{ 6}\phi^{ijk}y_i y_j y_k+
\fract{1}{ 24} \phi^{ijkl} y_i y_j y_k y_l. \nonumber
\end{align}
By substituting (\ref{eq:thetaineta}) and (\ref{eq:hiny}) for those in 
(\ref{eq:expdensity}), $f(y;\eta)$ is expressed by $y$, $\eta$, and
the derivatives of $\phi(\eta)$ at $\eta=0$. Then, (\ref{eq:uv-density})
is obtained from
$ f(\hu,\hv; \lambda)  =  f(\eta(\hu,\hv); (0,\ldots,0,\lambda))
J(\hu,\hv)$, where $J(\hu,\hv)$ is the Jacobian shown below.

The following lemma gives the Jacobian for the change of variables, which
is analogous to that for the tube formula
 of Weyl \nocite{bib:Weyl:1939:OVT} (1939).
The proof is due to
Satoshi Kuriki (pers. comm.); see Lemma 2.1 of Kuriki and Takemura
\nocite{bib:Kuri:Take:2000:SET} (2000) for the normal case.
\begin{lemma}  \label{lemma:jacobian}
The Jacobian $J(u,v)={\partial \eta / \partial (u,v)}$ is expressed as
\begin{eqnarray*}
 &&  \log J(u,v) \\
& \approx & -\fract{1}{2} \phi^{cpp} u_c
+(2d^{aa}-\phi^{aap})v
-\Bigl\{
2\dabdab -2d^{ab}\phi^{abp} +\fract{1}{2}\gabpgabp
\Bigr\} v^2\\
& &+\Bigl\{
\fract{1}{2}d^{cd}\phi^{ppp} -\fract{1}{4}\phi^{cdpp}
+\fract{1}{4}\phi^{cpp}\phi^{dpp} +\fract{1}{2}\phi^{acp}\phi^{adp}
+2d^{ac}(d^{ad}-\phi^{adp})
\Bigr\}u_cu_d\\
& & +\Bigl\{
6e^{aac} +d^{aa}\phi^{cpp} +4d^{ac}\phi^{app} -\phi^{aacp}
+\phi^{aad}\phi^{cdp} \\
& & \qquad +\fract{1}{2}\phi^{aap}\phi^{cpp} -2d^{cd}\phi^{aad}
-(2d^{ad}-\phi^{adp}) \phi^{acd}
\Bigr\}u_cv.
\end{eqnarray*}
\end{lemma}
{\it Proof.} Considering the projection of $d\eta_i = B_i^a(u) du_a +
B_i^p(u) dv + v{(\partial B_i^p/ \partial u_a)} du_a$ to the tangent
vectors and the normal vector, we get
\begin{equation}
\left(
 \begin{array}{c}
  B_i^a(u) \\
  B_i^p(u)
 \end{array}
\right)
\phi^{ij}(\eta(u))
 d\eta_j \approx
\left(
\begin{array}{cc}
 \phi^{ac}(u)& 0\\
 0 & 1
\end{array}
\right)
\left(
\begin{array}{cc}
 \delta_{cb}+vc_c^b(u)&0 \\
 v c^b_p(u) & 1
\end{array}
\right)
\left(\begin{array}{c}
 du_b\\ dv
\end{array}\right), \label{eq:weyldeta}
\end{equation}
where $c_b^a(u)$ and $c^a_p(u)$ are defined by $ {(\partial B_i^p/
\partial u_a)} = c_b^a(u) B_i^b(u) + c^a_p(u) B_i^p(u)$.  We have used a
simplified notation in which $\phi^{ac}(u)$, say, in the partitioned
matrix implies the $(p-1)\times(p-1)$ matrix instead of the
$(a,c)$-element of the matrix. It follows from the definition of the
metric that the determinant of the $p\times p$ matrix on the left hand
side of \refeq{weyldeta} is \[ \det \left\{\left(
 \begin{array}{c}
  B_i^a(u) \\
  B_i^p(u)
 \end{array}
\right) \phi^{ij}(\eta(u)) \right\} = \left\{\det(\phi^{ab}(u))
\det\left( \phi^{ij}(\eta(u)) \right)\right\}^{1\over2}. \] Hence we
find from \refeq{weyldeta} that the logarithm of the Jacobian is
expressed as $\log J(u,v) = \fract{1}{2}\log\det(\phi^{ab}(u)) + \log
\det(\delta_{ab}+vc_b^a(u)) -\fract{1}{2}\log\det\left(
\phi^{ij}(\eta(u)) \right)$.  Using the inverse matrix of $\phi^{ab}(u)$
with components $ \phi_{ab}(u) \approx \delta_{ab} - \phi^{abc} u_c - \{
4d^{ac}d^{bd} -2d^{ac}\phi^{bdp} -2d^{bd}\phi^{acp} -d^{cd}\phi^{abp}
-\phi^{ace}\phi^{bde} +\fract{1}{2}\phi^{abcd} \}u_c u_d$, the elements
of $c_b^a(u)$ is given by $ c_b^a(u) = {(\partial B_i^p/ \partial
u_a)}\, \phi^{ij}(\eta(u)) B_j^c(u) \phi_{cb}(u)$,  and therefore the
proof completes by noting the formal expansions
$ \log \det (I_p+A) = {\rm tr}(A) -\fract{1}{2} {\rm tr}(A^2) +\fract{1}{3}
 {\rm tr}(A^3) +\cdots$, and $(I_p+A)^{-1} = I_p - A + A^2 - A^3 +
 \cdots$ for $p\times p$ matrix $A$.

\subsection{Proof of Lemma~\ref{lemma:ftau}.}

Let the density function of $X$ be $\exp(\theta_X^i x_i
-\psi_X(\theta_X)-h_X(x)),$ where
$\theta_X=(\theta_X^1,\ldots,\theta_X^p)$ is the natural parameter and
$\mu=\partial \psi_X/\partial \theta_X$ is the expectation
parameter. The cumulant function of the simple sum $S=X_1+\ldots+X_m$ is
$m\psi_X(\theta_X)$ and the expectation parameter is $\eta_S =
m\mu$. The potential function of $S$ is then
$\phi_S(\eta_S)=\max_{\theta_X} \{ \theta_X'\eta_S - m\psi_X(\theta_X)\}
= m \phi_X(\eta_S/m)$. Since $Y=(\sqrt{n}/m)S$ has the same potential as
$S$, we obtain the potential of $Y$ as $\phi_Y(\eta)=m
\phi_X(\eta/\sqrt{n})$ with expectation
$\eta=(\sqrt{n}/m)\eta_S=\sqrt{n} \mu$. For $m=n$, it follows from the
assumption that $n\phi_X(\eta/\sqrt{n})=\phi(\eta)$.  This implies
$\phi_Y(\eta)=(m/n)\phi(\eta)=\phi(\eta)/\tau^2$ for $m>0$.

\subsection{Proof of Lemma~\ref{lemma:scaling}.}

We work on a random vector $\tilde Y = Y/\tau$, to which
the argument of the previous section is easily applied. The potential
function is $\tilde\phi(\tilde\eta,\tau)=\phi(\eta,\tau)$ with
expectation parameter $\tilde\eta=\eta/\tau$. Noting the $k$-th
derivative is $\partial^k \tilde\phi(\tilde\eta,\tau)/\partial
\tilde\eta^k = \tau^{k-2}\partial^k \phi(\eta,1)/\partial \eta^k$, and
especially the metric $\partial^2 \tilde\phi(\tilde\eta,\tau)/\partial
\tilde\eta_i \partial \tilde\eta_j = \partial^2 \phi(\eta)/\partial
\eta_i \partial \eta_j$ is $\delta_{ij}$ at $\tilde\eta=0$, 
Lemma~\ref{lemma:uv-density} is applied immediately to $\tilde Y$ 
 by the replacement (\ref{eq:reptilde}).
The tube-coordinates $(u,v)$ is replaced by $(\tilde u,\tilde
v)\leftrightarrow \tilde \eta$, which should be given by going back to
the specification of $\eta(u)$ and $B^p(u)$ using
$\tilde\phi(\tilde\eta,\tau)$ instead of $\phi(\eta,1)$. Fortunately a
simple transformation (\ref{eq:reptildeuv})
is consistent with the specification, and thus the joint density of
$(\tilde u, \tilde v)$ is easily brought back to that of $(u,v)$.  The
expression of $\log f(u,v;\lambda,\tau)$, although not shown to save the
space, is then obtained from (\ref{eq:uv-density}) by these
replacements, and by adding the logarithm of the Jacobian to it.

\subsection{Proof of Lemma~\ref{lemma:zc-formula}.}

We first consider the case $\tau=1$ to derive the cumulants
$\kappa_1,\kappa_2,\ldots,$ of $\hat w$.
The joint density of $(\hat u,\hat w)$ is
$ f(u,w; \lambda) = f(u,v(u,w); \lambda)(\partial v/ \partial w)$,
where the Jacobian is
$
 {\partial v/ \partial w} = \exp\{
 -\sum_{r\ge1} (c_r + u_c b_r^c)r w^{r-1} 
 -\fract{1}{2} ( \sum_{r\ge1} c_r r w^{r-1} )^2
\}$.
By rearranging the formula, $f(u,w; \lambda)$ is expressed in the form
\begin{eqnarray} \label{eq:fuw-a}
&& (2\pi)^{-{p-1\over2}} \exp\left(-\fract{1}{2} u_a u_a\right) \\
 &&\qquad \times
\exp\left(
A + A^a u_a + A^{ab} u_a u_b + A^{abc} u_a u_b u_c 
+ A^{abcd} u_a u_b u_c u_d
\right) + O(n^{-3/2}),\nonumber
\end{eqnarray}
where the coefficients $A=O(1)$, $A^{a}=O(n^{-1/2})$,
 $A^{ab}=O(n^{-1/2})$, $A^{abc}=O(n^{-1/2})$, and $A^{abcd}=O(n^{-1})$
are functions of $w$ and $\lambda$.  
We obtain the marginal density of $w$ 
by integrating (\ref{eq:fuw-a})
with respect to $u$. The logarithm of it is
$ A + A^{aa} + 3A^{aabb} 
+\fract{1}{2}A^a A^a + A^{ab}A^{ab}
+3 A^a A^{abb} +\fract{9}{2} A^{aac} A^{bbc} + 3A^{abc}A^{abc}
 + O(n^{-3/2})$, which becomes
$\log f(w; \lambda)
\approx -\fract{1}{2} \log(2\pi) -\fract{1}{2}(w-\lambda)^2
 +\Bigl[
\dabdab  -\fract{1}{2}d^{aa}d^{bb}  - d^{ab}\phi^{abp}
+\fract{1}{2}\gabpgabp
 +\fract{1}{2}\gappgapp +\fract{1}{6}(\phi^{ppp})^2
-\fract{1}{4}\phi^{aapp} -\fract{1}{8}\phi^{pppp} 
-d^{aa}\lambda  + \dabdab\lambda^2
-\fract{1}{2}d^{aa}\phi^{ppp}\lambda^2
-\fract{1}{3}\phi^{ppp}\lambda^3
-\fract{1}{8}\phi^{pppp}\lambda^4
+\fract{1}{8}\gappgapp\lambda^4
\Bigr]+w\Bigl[
d^{aa}+\fract{1}{2}\phi^{ppp}
+\Bigl\{
-d^{ab}\phi^{abp} +\fract{1}{2}d^{aa}\phi^{ppp} 
+\fract{1}{2}\gabpgabp +\fract{3}{8}\gappgapp
-\fract{1}{4}\phi^{aapp} \Bigr\}\lambda
+\fract{1}{2}\phi^{ppp}\lambda^2
+\fract{1}{6}\phi^{pppp}\lambda^3
-\fract{1}{4}\gappgapp\lambda^3
\Bigr]+w^2 \Bigl[
-\dabdab +d^{ab}\phi^{abp} -\fract{1}{2}\gabpgabp
-\fract{1}{2}\gappgapp -\fract{1}{4}(\phi^{ppp})^2
+\fract{1}{4}\phi^{pppp} +\fract{1}{4}\phi^{aapp}
+\fract{1}{8} \gappgapp\lambda^2
\Bigr]
 -\fract{1}{6}\phi^{ppp}w^3   -\fract{1}{24}\phi^{pppp}w^4
 +(\sum_i c_i w^i)\Bigl\{
w-\lambda -d^{aa} -\fract{1}{2}\phi^{ppp}+\fract{1}{2}\phi^{ppp}w^2
-\fract{1}{2}\phi^{ppp}\lambda^2
-\fract{1}{6}\phi^{pppp}\lambda^3
\Bigr\}  -\fract{1}{2}(\sum_i c_i w^i)^2 
 -\sum_i i c_i w^{i-1} -\fract{1}{2} (\sum_i i c_i w^{i-1})^2$.
For a parameter $t$, consider the formal expansion 
\begin{eqnarray}\label{eq:fwl-b}
 e^{wt} f(w; \lambda) &\approx& (2\pi)^{-{1\over2}}
\exp\left(-\fract{1}{2} ( w-(\lambda+t))^2 \right)\\
&&\quad \times \exp\left(
B_0 + B_1 w + B_2 w^2 + B_3 w^3 + B_4 w^4 
\right),\nonumber
\end{eqnarray}
where the coefficients $B_r$ are functions of $t$ and $\lambda$. 
By integrating (\ref{eq:fwl-b}) with respect to $w$, and taking the
logarithm of it, we obtain
$B_0 + B_2 + 3B_4
+\fract{1}{2} B_1^2  + B_2^2
+3 B_1 B_3 +\fract{15}{2} B_3^2
+(\lambda+t)\left(
B_1 + 3 B_3 +2 B_1 B_2  + 12 B_2  B_3 \right)
+(\lambda+t)^2\left(
B_2 +6 B_4 +3 B_1 B_3  +2 B_2^2
+18 B_3^2 \right)
+(\lambda+t)^3\left(B_3 +6 B_2 B_3\right)
+(\lambda+t)^4\left(B_4+\fract{9}{2}B_3^2\right)$.
This is compared with 
$\log \int e^{wt} f(w; \lambda)\,dw = t \kappa_1 +
{t^2\over2}\kappa_2 + {t^3\over6} \kappa_3 + {t^4\over 24}\kappa_4 +
\cdots$ to obtain
$ \kappa_1 \approx
  d^{aa} +c_0 +c_2+\Bigl\{
1 -2\dabdab +d^{ab}\phi^{abp}
-\fract{1}{2}d^{aa}\phi^{ppp}
-\fract{1}{2}\gabpgabp
-\fract{5}{8}\gappgapp
 +\fract{1}{4}\phi^{aapp}
+c_1 +2c_0c_2 + c_2(2d^{aa}-\phi^{ppp})
+3c_3 +6c_2^2\Bigr\}\lambda
+c_2 \lambda^2  +(c_3 +2c_2^2)\lambda^3$,
$\kappa_2  \approx 
1 - 2\dabdab +2d^{ab}\phi^{abp} -d^{aa}\phi^{ppp}
-\gabpgabp - \gappgapp
+\fract{1}{2}\phi^{aapp} 
+2c_1 +2c_2(2d^{aa}-\phi^{ppp})
+6c_3 +4c_0c_2 +14c_2^2
 +( -\phi^{ppp} +4c_2)\lambda 
 + \Bigl\{
\fract{1}{4}\gappgapp +(\phi^{ppp})^2
-\fract{1}{2}\phi^{pppp}
-4c_2\phi^{ppp} +6c_3 +16c_2^2
\Bigr\}\lambda^2$,
$\kappa_3  \approx 
-\phi^{ppp} +6c_2
+\Bigl\{
3(\phi^{ppp})^2 -\phi^{pppp}
-18c_2\phi^{ppp} +18c_3 +60c_2^2
\Bigr\}\lambda$,
$\kappa_4 \approx
3(\phi^{ppp})^2 -\phi^{pppp} -24c_2\phi^{ppp} +24c_3 +96c_2^2$,
and $\kappa_r\approx0$ for $r\ge5$.

The expression of $z_c(\hat w;\lambda,1)$ is obtained by applying the
Cornish-Fisher series (Johnson and Kotz, p.~66, 1994)
\nocite{bib:John:Kotz:Bala:1994:CUD} to these cumulants of $\hat w$. The
expression of $z_c(\hat w;\lambda,\tau)$ is then obtained by the
replacements of Lemma~\ref{lemma:scaling}, and by replacing $c_r$ with
$\tilde c_r=c_r \tau^{r-1}$.

\subsection{Proof of Lemma~\ref{lemma:change-origin}.}

Let us consider the local coordinates $\Delta\eta$ at $\eta(\hat u)$ so
that $\eta(\hat u)$ becomes the origin; \[ \eta_i = \eta_i(\hat u) +
B_i^j(\hat u)\Delta \eta_j,\quad i=1,\ldots,p \] defines the
reparameterization $\eta \leftrightarrow \Delta \eta$ at each $\hat u$.
The surface is now represented in the $\Delta \eta$-coordinates as
$\Delta \eta_p \approx-\hat d^{ab} \Delta \eta_a \Delta \eta_b - \hat
e^{abc} \Delta \eta_a \Delta \eta_b \Delta \eta_c$ using the Taylor
coefficients at $\eta(\hat u)$; $\hat d^{ab} \approx d^{ab} + (
\fract{1}{2} d^{ab}\phi^{cpp} +3e^{abc}) \hat u_c$ and $\hat e^{abc}
\approx e^{abc}$.  The derivatives $\hat\phi^{ij}$, $\hat\phi^{ijk}$,
and $\hat\phi^{ijkl}$ of $\phi(\eta)$ at $\eta(\hat u)$ with respect to
$\Delta\eta$-coordinates are obtained by comparing both sides of \[
\phi^{kl}(\eta)
B_k^i(\hat u) B_l^j(\hat u) = \hat \phi^{ij} +\hat\phi^{ijk}\Delta\eta_k
+\fract{1}{2}\hat\phi^{ijkl}\Delta\eta_k\Delta\eta_l+ \cdots. \] The
second derivatives are $\hat\phi^{ab}=\phi^{ab}(\hat u)$ of
(\ref{eq:phiab-u}), $\hat \phi^{ap} = 0$, and $\hat \phi^{pp}=1$. The
third derivative to the normal direction, i.e., $\hat\phi^{ppp}=-6\hat a$,
is (\ref{eq:hat-phippp}). The expressions for the other
$\hat\phi^{ijk}=\phi^{ijk}+O(n^{-1})$ are shown in Shimodaira
(2003). The forth derivatives are $\hat\phi^{ijkl}\approx \phi^{ijkl}$.

Although we assumed that $\phi^{ij}=\delta_{ij}$ at the origin in
(\ref{eq:uv-density}), the metric changes to $\hat\phi^{ij}$ at
$\eta(\hat u)$. We have to multiply the matrix $(\hat\phi^{ij})^{-1/2}$
to the vector $\Delta \eta$ to bring back the metric to the identity
matrix. This is equivalent to replacing $\delta_{ab}$ in the summation
over the indices $a,b$ with $ \hat\phi_{ab} =\phi_{ab}(\hat u)$
given in the proof of Lemma~\ref{lemma:jacobian}.

Consequently, $d^{aa}=d^{ab}\delta_{ab}$ in (\ref{eq:z-lambda-tau}) 
is replaced with $\hat d^{ab} \hat\phi_{ab}=\hat d_1$ of
(\ref{eq:hat-daa}). $\phi^{ppp}$ in (\ref{eq:z-lambda-tau}) is simply
replaced with $\hat\phi^{ppp}$ of (\ref{eq:hat-phippp}).
$O(n^{-1})$ terms change only $O(n^{-3/2})$; for example, $d_2=d^{ab}
d^{ab} = d^{ab} \delta_{bc} d^{cd} \delta_{da}$ becomes 
$\hat d_2=\hat
d^{ab}\hat\phi_{bc}\hat d^{cd}\hat\phi_{da} \approx \dabdab$.

\subsection{Proof of Lemma~\ref{lemma:distzq}.}

$\hat z_q(y)$ is expressed as $\hat w$, since $\bar c_r$ for $\hat z_q(y)$ is
that for $\hat z_\infty(y)$ plus $q_r$. Then the coefficients $c_r$ for
$\hat z_q(y)$ are calculated as $c_0 = -d^{aa}
-\fract{1}{6}\phi^{ppp}+q_0$, $c_1 = \dabdab -d^{ab}\phi^{abp}
+\fract{1}{2}d^{aa}\phi^{ppp} +\fract{1}{2}\gabpgabp
+\fract{1}{2}\gappgapp +\fract{17}{72}(\phi^{ppp})^2
-\fract{1}{4}\phi^{aapp} -\fract{1}{8}\phi^{pppp}
+\fract{1}{3}\phi^{ppp} (q_2-q_0) + q_1 + 2d^{aa} q_2 -2q_0 q_2$, $c_2 =
\fract{1}{6}\phi^{ppp} + q_2$, $c_3 = -\fract{5}{72} (\phi^{ppp})^2
+\fract{1}{24} \phi^{pppp} -\fract{2}{3}\phi^{ppp} q_2 - 2 q_2^2 + q_3$.
  By applying Lemma~\ref{lemma:zc-formula}
to $z_c(\hat w;\lambda,1)$ with these coefficients, we obtain the
distribution function of $\hat z_q(y)$.

\subsection{Proof of Lemma~\ref{lemma:threestepboot}.}

We derive the expression of $\tilde
z_2(\eta(0,\lambda),\tau_1,\tau_2)$ so that  $\tilde
z_2(y,\tau_1,\tau_2)$ is obtained by the replacements of
Lemma~\ref{lemma:change-origin}. By letting $y=\eta(0,\lambda)$, the
right hand side of (\ref{eq:z2-integration}) becomes
\begin{align}\label{eq:znbintinv} 
\tilde z_2(\eta(0,\lambda),\tau_1,\tau_2)=
\Phi^{-1}\Bigl\{
\int\!\!\int &\Phi(\tilde z_1(\eta(\hu,\hv),\tau_2))
 f(\hu,\hv;\lambda,\tau_1)\,d\hu\,d\hv
\Bigr\} \\ &\approx \Phi^{-1}\Bigl\{
\int \Phi(\tilde z_1(\eta(0,\hv),\tau_2)) f(\hv;\lambda,\tau_1)\,d\hv
\Bigr\},  \nonumber
\end{align}
where the last approximation
follows from the fact that $\hu$ and $\hv$
 are approximately independent ignoring $O(n^{-1/2})$ terms;
$\log f(\hu,\hv;\lambda,\tau) = -\fract{1}{2}p\log(2\pi\tau)
-\fract{1}{2}(\hv-\lambda)^2\tau^{-2} -\fract{1}{2}(\hu_a)^2\tau^{-2}
 + O(n^{-1/2})$,
and
hence terms of $\bar b_r^c=O(n^{-1})$ in \refeq{winv} for $\tau_2 \tilde
z_1(y,\tau_2)$ contribute only $O(n^{-3/2})$ in the integration.

To carry out the integration of \refeq{znbintinv}, we consider a
transformation described below.  Let us define a generalization of the
pivot by $z_\infty(\hat u,\hat v;\lambda,\tau) := z_c(\hat v;\hat
u,\lambda,\tau)$ with all $c_r=0$. This reduces to $\hat z_\infty(y)$
 when $\lambda=0$ and $\tau=1$, and a similar argument as
\refsec{pivot} shows that
$\Pr\{ z_\infty(\hat U,\hat V;\lambda,\tau) \le x;\lambda,\tau \}
\approx \Phi(x)$.
By solving $z_\infty(\hat u,v;\lambda,\tau)=z$ for $v$, we define the
inverse function  $v_\infty(\hat u,z;\lambda,\tau)$, which satisfies
$z_\infty(\hat u,v_\infty(\hat u,z;\lambda,\tau);\lambda,\tau)=z$. 
Using \refeq{vinw}, we may obtain
$v_\infty(0,z;\lambda,\tau)\approx \lambda  + \tau z \biggl[
1 - \fract{1}{2}\phi^{ppp}\lambda + 
(\fract{1}{8}\gappgapp + \fract{3}{8}(\phi^{ppp})^2 
- \fract{1}{4}\phi^{pppp})\lambda^2
\biggr]
+\tau^2 \biggl[
d^{aa} + \fract{1}{6}\phi^{ppp}
+\Bigl( 
-2\dabdab +\fract{1}{4}\phi^{aapp} + d^{ab}\phi^{abp}
-\fract{1}{2}\gabpgabp -\fract{5}{8}\gappgapp
 -\fract{1}{2} d^{aa}\phi^{ppp} -\fract{1}{3}(\phi^{ppp})^2
+\fract{1}{6}\phi^{pppp}
\Bigr)\lambda + z^2 \Bigl(
-\fract{1}{6}\phi^{ppp} 
+(\fract{1}{3}(\phi^{ppp})^2 - \fract{1}{6}\phi^{pppp})\lambda
\Bigr)
\biggr]
+\tau^3\biggl[
z\Bigl(
-\dabdab +\fract{1}{4}\phi^{aapp} +d^{ab}\phi^{abp}
-\fract{1}{2}\gabpgabp -\fract{1}{2}\gappgapp
-\fract{1}{2}d^{aa}\phi^{ppp}
 -\fract{17}{72}(\phi^{ppp})^2
+\fract{1}{8}\phi^{pppp}
\Bigr)
+z^3(\fract{5}{72}(\phi^{ppp})^2 -\fract{1}{24}\phi^{pppp})
\biggr]$
.

Considering $\tilde z_1(\eta(0, v_\infty(0,z;\lambda,\tau_1)),\tau_2) =
(\lambda + z \tau_1) \tau_2^{-1} + O(n^{-1/2})$, and $\Phi(x +\delta)
\approx \Phi(x) + f(x) (\delta - \fract{1}{2}x \delta^2)$ for $x=O(1)$
and $\delta=O(n^{-1/2})$, we find that the integration \refeq{znbintinv}
is now expressed as
\begin{align} \label{eq:znbintinz}
\Phi^{-1}\Bigl\{
\int &\Phi\Bigl(\tilde z_1(\eta(0,
v_\infty(0,z;\lambda,\tau_1)),\tau_2)\Bigr) f(z)\,dz
\Bigr\}\\
&\approx \Phi^{-1}\Bigl\{
\int \Bigl(
 \Phi(az+b) + f(az+b)
\sum_{r=0}^5 A_r z^r
\Bigr) f(z)\,dz
\Bigr\},\nonumber
\end{align}
where $f(z)$ is the standard normal density function, and $a=\tau_1
\tau_2^{-1}, b=\lambda\tau_2^{-1}$, and $A_r=O(n^{-1/2})$ are
independent of $z$.  The expressions for $A_r$ are not shown here to
save the space, but they are functions of the geometric quantities as
well as $\lambda$, $\tau_1$, and $\tau_2$.  Let $c=b(1+a^2)^{-1/2}$.
Using
$\int_{-\infty}^{\infty} \Phi(az+b) f(z)\,dz = \Phi(c)$
and $g_r=\int_{-\infty}^{\infty} z^r f(az+b) f(z)\,dz/f(c)$,
 the integration \refeq{znbintinz} is expressed as
\begin{equation} \label{eq:znbinting}
\Phi^{-1}\Bigl\{
\Phi(c) + f(c)\sum_{r=0}^5 A_r g_r
\Bigr\},
\end{equation}
where $g_0=(1+a^2)^{-1/2}$,
$g_1= -ab(1+a^2)^{-1}g_0$, $g_2=(1+a^2(1+b^2))(1+a^2)^{-2}g_0$,
$g_3=-ab(3+a^2(3+b^2))(1+a^2)^{-3}g_0$,
$g_4=(3+6a^2(1+b^2)+a^4(3+6b^2+b^4))(1+a^2)^{-4}g_0$, and
$g_5=-ab(15+10a^2(3+b^2)+a^4(15+10b^2+b^4))(1+a^2)^{-5}g_0$.
To derive $g_r$, we may use
$\int_{-\infty}^{\infty} z^{2r} f(z)\,dz = (2r)!/(2^rr!)$,
which becomes $1,3,15,105,945$, say,  for $r=1,2,3,4,5$.

By noting $\Phi^{-1}(\Phi(x) + f(x) \delta ) \approx x + \delta +
\fract{1}{2} x \delta^2$ for $x=O(1)$ and $\delta=O(n^{-1/2})$, we
finally obtain from \refeq{znbinting} that
\begin{equation} \label{eq:znbintinab}
 \tilde z_2(\eta(0,\lambda),\tau_1,\tau_2)\approx
c\Bigl\{
1+ \fract{1}{2}\Bigl(\sum_{r=0}^5 A_r g_r \Bigr)^2
\Bigr\} + \sum_{r=0}^5 A_r g_r.
\end{equation}
It is straightforward, though very tedious, to verify that
\refeq{znbintinab} is in fact equivalent to
$\zeta_3(\gamma_1,\ldots,\gamma_6,\tau_1,\tau_2,0)$ ignoring
$O(n^{-3/2})$ terms.

The proof for the three-step bootstrap goes similarly as that for the
two-step bootstrap. First we note that $\tilde
z_3(\eta(0,\lambda),\tau_1,\tau_2,\tau_3) \approx \Phi^{-1}\Bigl\{ \int
\Phi(\tilde z_2(\eta(0, v_\infty(0,z;\lambda,\tau_1)),\tau_2,\tau_3))
f(z)\,dz \Bigr\} $, and that $ \tilde z_2(\eta(0,
v_\infty(0,z;\lambda,\tau_1)) ,\tau_2,\tau_3) = (\lambda + z \tau_1)
(\tau_2^2 + \tau_3^2)^{-1/2} + O(n^{-1/2})$.  Then $\tilde
z_3(\eta(0,\lambda), \tau_1,\tau_2,\tau_3)$ is expressed in the same
form as \refeq{znbintinab} but using $ a=\tau_1 (\tau_2^2 +
\tau_3^2)^{-1/2}$, $b=\lambda (\tau_2^2 + \tau_3^2)^{-1/2}$, and
different $A_r$'s. A straightforward calculation shows that this is
equivalent to $\zeta_3(\gamma_1,\ldots,\gamma_6,\tau_1,\tau_2,\tau_3)$
ignoring $O(n^{-3/2})$ terms.

\section{Mathematrica session.}

\setcounter{section}{0}
\renewcommand\thesection{*\arabic{section}}

\newcommand{\m}[1]{{\tt #1}}
\newcommand{\e}[1]{{\rm #1}}
\newcommand{\oi}{O(n^{-1/2})}
\newcommand{\oii}{O(n^{-1})}
\newcommand{\oiii}{O(n^{-3/2})}

\newcommand{\MN}{{\it Program}}
\newcommand{\SH}{{\it Paper}}

This section explains briefly Shimodaira (2003), which is a computer
program in {\it Mathematica} notebook document referred to as \MN\ here.
This program is available in the {\it Mathematica} notebook format, PDF,
and HTML from the author. \MN\ proves the third-order accuracy of the
bias-corrected $p$-value calculated by the three-step multiscale
bootstrap resampling.  This newly devised bootstrap method is described
in the paper Shimodaira (2004) referred to as \SH\ here.

The document starting below has the same section structure as \MN\ with
a brief explanation. This helps us to find the appropriate results in
\MN.  The section numbers of \MN\ are indicated with asterisk~(*) to
avoid confusion.  Notational differences between \SH\ and \MN\ are also
explained.

\MN\  consists of three parts, ``Exponential Family of Distributions,''
``Tube-Coordinates and $z_c$-formula,'' and ``Bootstrap Methods.''  Each
part is an independent {\it Mathematica} session, and should be run
separately. In the first two parts, the tensor notation is heavily used.
The add-on package {\it MathTensor} is required to run the session by
yourself. I have first hand-calculated the results of the first two parts
involving tensor notation, and later used {\it MathTensor} to verify the
results.  The last part has been calculated only by {\it Mathematica}.


\section*{*Part I\\  Exponential Family of Distributions}

\section{Startup}
This section initializes the {\it Mathematica} session.
\subsubsection{packages}

\m{Statistics`ContinuousDistributions`} and \m{MathTensor} are loaded.

\subsubsection{error messages}
\subsubsection{distribution functions}

\m{f[x]}, \m{F[x]}, and \m{Q[p]} are the density function, the
distribution function, and the quantile function, respectively, for the
standard normal random variable. $\m{F[x]}$ and $\m{Q[p]}$ are denoted
$\Phi(x)$ and $\Phi^{-1}(p)$, respectively, in \SH. \m{f[x]} is also
denoted as $f(x)=(2\pi)^{-1/2}\exp(-x^2/2)$ here.

\section{Normal distribution}

\subsection{The moments of the multivariate normal distribution}

\subsubsection{the central moments of the standard normal variable in
   one dimension}
\[
\m{intx2f[n]}= \int_{-\infty}^{\infty} x^{2n} f(x)\,dx=\frac{(2n)!}{2^n n!}
\]
This gives 1, 3, 15, 105, 945, \ldots, for $n=$ 1, 2, 3, 4, 5, \ldots.

\subsubsection{define tensors}
\subsubsection{the central moments for the multivariate case}

Let $x=(x_a)$ denote a multivariate normal random vector with mean 0 and
covariance identity. Then, the central moments, $\alpha_{ab}=E(x_a
x_b)$, $\alpha_{abcd}=E(x_a x_b x_c x_d)$, and $\alpha_{abcdef}=E(x_a
x_b x_c x_d x_e x_f)$ are given by \m{kd2[a,b]}, \m{kd4[a,b,c,d]}, and
\m{kd6[a,b,c,d,e,f]}, respectively. They satisfy, for example,
$\alpha_{112222}=E(x_1^2 x_2^4)=\m{intx2f[1] intx2f[2]}=3$.

\subsection{The expectation of the exponential of polynomial functions} 
\label{sec:logeexppoly}

\subsubsection{central case}

Let $\e{poly}(x)$ be a polynomial of $x$, 
\[
 \e{poly}(x) \approx
 a0 +  a1^a x_a +  a2^{ab} x_a x_b +  a3^{abc} x_a x_b x_c
 +  a4^{abcd} x_a x_b x_c x_d,
\]
where $a0$ is $O(1)$, $a1, a2, a3$ are $O(n^{-1/2})$, and $a4$ is
$O(n^{-1})$.  We use '$\approx$' to indicate equivalence up to $\oii$
terms ignoring the error of $\oiii$.  The polynomial is expressed in
\MN\ as
\[
 \e{poly}(x) =
 a0 +o \left( a1^a x_a +  a2^{ab} x_a x_b +  a3^{abc} x_a x_b x_c\right)
 + o^2 a4^{abcd} x_a x_b x_c x_d,
\]
where the $\m{o}$ indicates a $\oi$ term, and $\m{o^2}$ indicates
a $\oii$ term. Then,
\[
\log E\left\{\exp(\e{poly}(x))\right\} \approx \m{logeexppolyb0}.
\]
\subsubsection{noncentral case}

When $E(x_a)=b_a$, we obtain
\[
\log E\left\{\exp(\e{poly}(x))\right\} \approx \m{logeexppoly}.
\]

\section{Exponential family}

\subsection{The standard form}

\subsubsection{simplification functions}

\subsubsection{define tensors}

\subsubsection{the log of the density function}

Let $y=(y_a)$ be a multivariate random variable of dimension \m{dim},
and the density function is specified by \[ f(y;\theta) = \exp(\theta^a
y_a - h(y) - \psi(\theta)), \] which corresponds to
eq.~\refeq{expdensity} of \SH\ given below \[ \exp(\theta^i y_i -
h(y) - \psi(\theta)).  \] The expression is given in $\log
f(y;\theta)=\m{logdensity}$ in \MN.  The dimension is denoted $p$ in
\SH, but it is denoted either \m{dim} or \m{9} in \MN; see
Table~\ref{tab:index}.  So $y=(y_1,\ldots,y_9)$,
$\theta=(\theta^1,\ldots,\theta^9)$ here.  We made \m{dim}=\m{9} only
for a technical reason of {\it MathTensor}, and our calculation is
independent of this choice.  The indices to run $1,...,p$ in \SH\ are
$i,j,\ldots$, but those to run $1,\ldots,9$ in \MN\ are
$\m{a},\m{b},\ldots$ for free indices, and $\m{p}, \m{q}, \ldots$ for
dummy indices.  It should be noted that the indices to run $1,...,p-1$
in \SH\ are $a,b,\ldots$, but those to run $1,...,8$ in \MN\ are
$\m{a'}, \m{b'}, \ldots$ or $\m{p'},\m{q'},\ldots$.  These differences
are quite confusing, but a restriction using {\it MathTensor}.

In the tensor notation, covariant and contravariant indices are
distinguished. In other words, a subscript and a superscript are summed
over their range if they are a matched pair of indices. {\it MathTensor}
and \MN\ follow this rule. In \SH, however, the summation convention
takes place for a pair of the same indices irrespective of subscripts
and superscripts.  This does not give any difference in the calculation,
because we formally set the metric matrix to be identity, and specify
the metric explicitly by the matrix multiplication when necessary.

\begin{table}[htbp]
 \caption{Notations of the indices} \label{tab:index}
\begin{center}
 \begin{tabular}{ccc}
\hline
 \SH & \MN & description\\
\hline
$p$  & \m{dim} or 9 & dimensions of the parameter vector\\
$i,j,\ldots$ & $a,b,\ldots$ and $p,q,\ldots$& indices run over
  $1,\ldots,p$ (=9)\\
$a,b,\ldots$ & $a',b',\ldots$ and $p',q',\ldots$ & indices run over
$1,\ldots,p-1$ (=8) \\
\hline
 \end{tabular}
\end{center}
\end{table}

The natural parameter vector is $\theta=(\theta^a)$, and
the expectation parameter vector $\eta=(\eta_a)$ is defined by
$\eta_a=E(y_a;\theta)$.
The potential function $\phi(\eta)$ is defined by
$\phi(\eta)=\max_{\theta}\{\theta^a \eta_a - \psi(\theta)\}$, and the
two parameterization are related to each other by $\eta_a={\partial
\psi/\partial \theta^a}$ and $\theta^a = {\partial \phi/ \partial
\eta_a}$. Without loss of generality, we assume
\[
 \frac{\partial \phi}{\partial\eta_a}\biggr|_0 = 0,\quad
 \frac{\partial^2 \phi}{\partial\eta_a\partial\eta_b}\biggr|_0 =
 \delta^{ab}=\m{Kdelta[a,b]}.
\]
For higher order derivatives of $\phi$ at $\eta=0$, we denote
\[
 \phi^{abc} = 
 \frac{\partial^3
 \phi}{\partial\eta_a\partial\eta_b\partial\eta_c}\biggr|_0=o \m{\phi3^{abc}} ,
 \quad
 \phi^{abcd} = 
 \frac{\partial^4
 \phi}{\partial\eta_a\partial\eta_b\partial\eta_c\partial\eta_d}\biggr|_0
= o^2 \m{\phi4^{abcd}},
\]
where the $\m{o}$ indicates a $\oi$ term, and $\m{o^2}$ indicates
a $\oii$ term.

\subsection{The expression of $\psi(\theta)$ in terms of $\eta$}

\subsubsection{derivation}
\subsubsection{result}

\[
 \psi(\theta) \approx \psi(0) + \frac{1}{2} \eta_p \eta_p + \frac{1}{3}
\eta_p \eta_q \eta_r \phi^{pqr} + 
\frac{1}{8} \eta_p \eta_q \eta_r \eta_s\phi^{pqrs} = \m{psieta},
\]
which corresponds to eq.~\refeq{thetaineta}. The indices \m{p},
\m{q},\ldots are dummy indices to run $1,...,9$ in \MN.

\subsection{The expression of $h(y)$ in terms of $\phi$ derivatives}

\subsubsection{derivation}
\subsubsection{result}

\begin{align*}
h(y)\approx &\frac{1}{2} \m{dim} \log(2\pi) - \psi(0) +\frac{1}{2}y_p y_p
 - \frac{1}{2} y_p \phi^{pqq}+
 \frac{1}{6} y_p y_q y_r \phi^{pqr} - \frac{1}{6}(\phi^{pqr})^2 +\\
&\frac{1}{4} y_p y_q \phi^{prs}\phi^{qrs} + \frac{1}{8}\phi^{ppqq} -
 \frac{1}{4}y_p y_q \phi^{pqrr} + \frac{1}{24}y_p y_q y_r y_s
 \phi^{pqrs}
=\m{hinyphi},
\end{align*}
which corresponds to eq.~\refeq{hiny}.

\subsection{The canonical form}

\label{sec:logdensityy}

\subsubsection{the summary of the previous sections}
\subsubsection{result}

Using the above expressions for $\psi(\theta)$ and $h(y)$, we obtain the
expression for the density function using $\eta$-parameter.
\[
\log f(y;\eta) = \m{logdensityy}.
\]
The metric, i.e., the second derivative of $\phi$, evaluated at a
general point $\eta$ is given by
\[
 \frac{\partial^2 \phi}{\partial \eta_a \partial \eta_b} = \m{phi2eta}.
\]

\section*{*Part II\\ Tube-Coordinates and $z_c$-formula}

\section{Startup}
This section initializes the {\itshape Mathematica} session.

\subsubsection{packages}
\subsubsection{error messages}
\subsubsection{distribution functions}

\section{Exponential family}

\subsection{The expectation of the exponential of a polynomial function of the normal vector}
\subsubsection{define tensors}
\subsubsection{logeexppoly}

This is a result of Section~\ref{sec:logeexppoly}.

\subsection{The canonical form of the density function}
\subsubsection{define tensors}

\subsubsection{logdensityy}

This is a result of Section~\ref{sec:logdensityy}.

\subsubsection{phi2eta}

This is a result of Section~\ref{sec:logdensityy}.

\section{Tube-coordinates}
\subsection{Preliminary}

\subsubsection{indices}
\subsubsection{simplification functions} \label{sec:simpfunc}
\subsubsection{define tensors}

\subsection{The coordinates around the smooth surface}
\subsubsection{smooth surface}

Let $u=(u_{a'})=(u_1,\ldots,u_8)$ parameterize a surface in 9-dimensional
space. The surface $\partial \cR=\{\eta(u)\}$ is specified by
\[
 \eta_{a'}(u)=u_{a'},\quad
 \eta_9(u) \approx - o d^{a'b'}u_{a'}u_{b'} - o^2 e^{a'b'c'}u_{a'}u_{b'}u_{c'},
\]
where the curvature matrix $o d^{a'b'}$ is $\oi$, and $o^2 e^{a'b'c'}$
is $\oii$.

\subsubsection{tangent vectors}

The $b$-th element of the $a'$-th tangent vector
of $\partial\cR$ at $\eta(u)$ is
\[
 B_{b}^{a'} = \frac{\partial \eta_b(u)}{\partial u_{a'}}
\]
for $a'=1,\ldots,8$, and $b=1,\ldots,9$. The expressions are given in
\m{foo3} and \m{foo4}, respectively, for $b=b'$ and $b=9$. They are
obtained by using the differential operator \m{difa} defined in
Section~\ref{sec:simpfunc}.
The metric for the tangent vectors is
\[
 \phi^{a'b'}(u)=\frac{\partial^2\phi(\eta)}{\partial \eta_p \partial
 \eta_q}\biggr|_{\eta(u)} B_p^{a'}(u) B_q^{b'}(u) =\m{phi2bu}.
\]

\subsubsection{the normal vector}

The $a$-th element of the normal vector at $\eta(u)$ is $ B_{a}^{9}$ for
$a=1,\ldots,9$. The expressions are given in \m{foo15} and \m{foo16},
respectively, for $a=a'$ and $a=9$. $B_a^9$ is orthogonal to the tangent
vectors with respect to the metric defined by the second derivative of
$\phi(\eta)$, and it has the unit length directing away from the region
$\cR$.

\subsubsection{$(u,v)$-coordinate system}

We define the reparameterization $\eta \leftrightarrow (u,v)$ by
\[
 \eta_a(u,v) = \eta_a(u) + B_a^9(u)v,\quad
a=1,\ldots,9,
\]
where $v$ is a scalar parameter. The expressions of $\eta_a(u,v)$ are
given in \m{foo21} and \m{foo22}, respectively, for $a=a'$ and $a=9$.

\subsection{Change of variables}
\subsubsection{Jacobian}

The Jacobian of the change of variables $\eta\leftrightarrow (u,v)$ is
\[
 J=\det\left(\frac{\partial \eta(u,v)}{\partial (u,v)}\right).
\]
The asymptotic expression of $\log J$ in terms of $(u,v)$ is given in
\m{logdetJ}.  This has been obtained in \reflemma{jacobian} of
\refsec{lemma1} using a sophisticated argument, but here \m{logdetJ} is
obtained more directly from the definition using {\it MathTensor}.

\subsubsection{density function $f(u,v|v0)$}

The random vector $y$ is also represented in $(u,v)$-coordinates. The
density function $f(u,v|v0)$ is expressed by $\log f(u,v|v0)=
\m{logdensityuv}$ with the true parameter vector $\eta=(0,\ldots,v0)$.
This corresponds to eq.~\refeq{uv-density} of \SH.  The scalar parameter
$v0$ is denoted $\lambda$ in \reflemma{uv-density} of \SH.

\section{$z_c$-formula}

\subsection{Modified signed distance $w$}

\subsubsection{define the modified signed distance as a series of $v$}

Let us define a modified signed distance $w$ by
\[
 w = v+\sum_{r=0}^\infty \e{cbr}[r] v^r + u_{a'}\sum_{r=0}^\infty
 \e{br}[a',r] v^r,
\]
where the coefficients $\e{cbr}[r]$ and $\e{br}[a',r]$ are denoted by
$\bar c_r$ and $\bar b_r^c$, respectively, in eq.~\refeq{winv} of \SH.
We assume that cbr[0] and cbr[2] are $\oi$, and cbr[1], cbr[3], and all
$\e{br}[a',r]$ are $\oii$. The other cbr[r] with $r\ge4$ are all
$\oiii$.

The inverse series specifies
\[
 v = w - \sum_{r=0}^\infty \e{cr}[r] w^r - 
u_{a'}\sum_{r=0}^\infty  \e{br}[a',r] w^r = \m{vinuw},
\]
where $\e{cr}[r]$ has the same order as $\e{cbr}[r]$.  Taking into
account the magnitude of the coefficients, the expression is actually
given in
\[
 \m{vinuw} = w - o(\e{cr}[0]+w^2 \e{cr}[2]) -
 o^2(\e{cr}[1]+w^3\e{cr}[3] + u_{a'} \e{br}^{a'})
\]
in \MN, where $\e{br}^{a'}$ indicates the $\oii$ polynomial functions of
$w$. The relation between these two sets of coefficients are given in \m{rule47}.

\subsubsection{density function of $w$}

The joint density function of $(u,w)$ is $f(u,w|v0)$. The asymptotic
expression is given in $\log f(u,w|v0)=\m{logdensityuw}$. The marginal
density $f(w|v0)=\int f(u,w|v0)\,du$ is obtained by using
\m{logeexppoly}, and the expression is given in $\log
f(w|v0)=\m{logdensityw}$. These expressions are found in the proof of
\reflemma{zc-formula} of \SH.

\subsubsection{cumulants of $w$}

By applying \m{logeexppoly} to \m{logdensityw}, we calculate
$\log \int e^{wt} f(w; \lambda)\,dw = t \kappa_1 +
{t^2\over2}\kappa_2 + {t^3\over6} \kappa_3 + {t^4\over 24}\kappa_4 +
\cdots$  to obtain
$\m{cumulantw}=\{\kappa_1,\kappa_2,\kappa_3,\kappa_4\}$, the cumulants
of $w$.

\subsection{Distribution function of $w$}

\subsubsection{Cornish-Fisher expansion  (p.66 of Johnson and Kotz
   1994)}

The Cornish-Fisher expansion for the standardized random variable is
taken from Johnson and Kotz (1994) as shown in \m{cfexpx},
and  the same expansion for nonstandardized variable is obtained as
shown in \m{cfexpw}.

\subsubsection{Cornish-Fisher expansion of $w$}

By applying \m{cfexpw} to \m{cumulantw}, we obtain
\[
 \m{zformula} = \Phi^{-1}(\Pr\{W\le w;v0\}).
\]
This corresponds to the $z_c$-formula $z_c(\hat w;\lambda,\tau)$ with
$\hat w=w$, $\lambda=v0$, and $\tau=1$ obtained in \reflemma{zc-formula}
of \SH.

\subsubsection{$z_c$-formula using a simplified notation}

The same expression as \m{zformula}, but without {\it MathTensor}
 notation, is given in \m{zform}. The following table shows the notational
 differences.
\begin{table}[h]
\caption{Notations of the geometric quantities} \label{tab:geometric}
\begin{center}
\begin{tabular}{ccc}
\hline
\SH & \m{zformula} & \m{zform}\\
\hline
$d^{aa}$ &  ${d_{p'}}^{p'}$ & Daa\\
$(d^{ab})^2$ & $({d_{p'}}^{q'}) ({d_{q'}}^{p'})$ & Dab2\\
$\phi^{ppp}$ & $\phi3^{999}$ & P999\\
$\phi^{pppp}$ & $\phi4^{9999}$ & P9999\\
$(\phi^{app})^2$ & $({\phi3^{99}}_{p'}) (\phi3^{99p'})$ & P99a2\\
$d^{ab}\phi^{abp} $& $ (d^{p'q'}) ({\phi3^9}_{p'q'})$ &  DabP9ab\\
 $(\phi^{abp})^2$ & $ ({{\phi3^9}_{p'}}^{q'})({{{\phi3}^9}_{q'}}^{p'})$ &
 P9ab2\\
$\phi^{aapp}$  & ${{\phi4^{99}}_{p'}}^{p'}$ & P99aa\\
\hline
\end{tabular} 
\end{center}
\end{table}

\subsubsection{scaling by the factor $\tau$}

By applying the scaling rule described in \reflemma{scaling} of \SH, we
obtain the expression of $z_c(w;v0,\tau)$ in \m{zformulatau} or
\m{zformtau}.

\subsection{Local coordinates at the projection}

\subsubsection{the expression for the surface in the local coordinates}

We consider a local coordinate $\Delta\eta=(\Delta \eta_a)$ around a
point $\eta(u0,0)$ on the surface, where $u0$ indicates any specified
value of $u$. We will use a particular value of $u0$ specifying the
projection of $y$ onto the surface. This value is denoted $\hat u$ in
\reflemma{change-origin} of \SH. The change of variables
$\eta\leftrightarrow \Delta\eta$ is specified by
\[
 \eta_a = \eta_a(u0,0) + B_a^b(u0) \Delta \eta_b,\quad
  a=1,\ldots,9.
\]
The expression for $\eta_{a'}$ is given in \m{foo93} and that for
 $\eta_9$ is in \m{foo94}. The surface $\partial\cR$ is now expressed as
$\Delta\eta_9 = -\hat d^{a'b'} \Delta\eta_{a'}\Delta\eta_{b'}
 -\hat e^{a'b'c'} \Delta\eta_{a'}\Delta\eta_{b'} \Delta\eta_{c'}$, where
 the expression for $\hat d^{a'b'}$ is in \m{foo101} and that for $\hat
 e^{a'b'c'}$ is in \m{foo102}. These expressions are also found in the
 proof of \reflemma{change-origin}.

\subsubsection{the expressions for the potential derivatives}

The expression is given for
 \[
{\partial^2 \phi(\eta)\over \partial \eta_p \partial \eta_q
}\biggr|_{\eta=\eta(u0)+ B(u0)\Delta\eta}
B_p^a(u0) B_q^b(u0) = \m{foo114[a,b]}. \]

\subsubsection{Geometric quantities at the projection}

By comparing both sides of
 \[
\m{foo114[a,b]}
= \hat \phi^{ab} +\hat\phi^{abc}\Delta\eta_c
+\fract{1}{2}\hat\phi^{abcd}\Delta\eta_c\Delta\eta_d+ \cdots, \]
we obtain the expressions for $\hat\phi^{ab}$, $\hat\phi^{abc}$, and
 $\hat\phi^{abcd}$ as shown in \m{foo121}.
The expression for the inverse matrix of the metric
$\hat\phi_{a'b'}=(\hat\phi^{a'b'})^{-1}$ is shown in \m{foo131}, and it
 is used for $\hat\phi_{a'b'}\hat d^{a'b'}=\m{foo132}$.

\subsubsection{$z_c$-formula}

The $z_c$-formula evaluated by taking $\eta(u0,0)$ as the origin of the
local coordinates is shown in \m{zformulau0} or \m{zformulatauu0}. They
correspond to $z_c(w;u0,v0,1)$ and $z_c(w;u0,v0,tau)$ in the notation of
\reflemma{change-origin} of \SH. From the expressions, we observe that
they depend on $u0$ only linearly by ignoring $\oiii$ terms, and the
magnitude is $\oii$.

\section*{*Part III\\ Bootstrap Methods}

\section{Startup}
This section initializes the {\itshape Mathematica} session.

\subsubsection{packages}

\subsubsection{error messages}

\subsubsection{distribution functions}

\section{Asymptotic Analysis of Bootstrap Methods}

This section calculates the distribution functions of $z$-values
appearing in several bootstrap methods for showing their asymptotic
accuracies.

\subsection{Preliminary}

\subsubsection{simplification functions}

\subsubsection{zc-formula}

\m{zc[w,\{c0,c1,c2,c3\},v0,tau]=zformtau} is
$z_c(w;v0,tau)=\Phi^{-1}(\Pr\{ W\le w; v0,tau\})$, where $u0$ is assumed
zero. The coefficients $\{c0,c1,c2,c3\}$ specify the modified signed
distance $w$, where $c0=\e{cr}[0],\ldots, c3=\e{cr}[3]$ in the notation of
\m{zformtau}.
We do not need to use the more general \m{zformulatauu0}
in which $u0\neq0$, because the $u0$ terms contribute only $\oiii$ in
the results of our calculation below.

\subsubsection{modified signed distance}

Here we write $cb0=\e{cbr}[0],\ldots,cb3=\e{cbr}[3]$. There are
relations between $\{cb0,\ldots,cb3\}$ and $\{c0,\ldots,c3\}$, and
\m{cb2cc} calculates the latter from the former.

\subsection{The pivot and some existing bootstrap methods}

\subsubsection{pivot statistic}

The pivot is defined as \[ \hat z_\infty(y)=-\Phi^{-1}(\Pr\{ V\ge \hat
v; u0=\hat u,v0=0\}) \] in \refsec{pivot} of \SH, where $(\hat u,\hat
v)$ is the tube-coordinates of $y$.  This is denoted as
\m{z8[v]}=\m{zc[v,\{0,0,0,0\},0,1]} in \MN\ for $\hat u=u0=0$, $\hat
v=\m{v}$.  The $\hat z_q(y)$ in \reflemma{distzq} of \SH\ is denoted
$\m{zq}[v] = \m{z8}[v] + o(q0 + q2 v^2 )+ o^2 (q1 v + q3 v^3
)$. \m{zq2qq} calculates \m{qq}=\{q0,q1,q2,q3\} from any
$z$-value. $\m{cb}_r$'s and $\m{c}_r$'s for \m{z8} are in \m{cbzq} and
\m{czqq}, respectively. The distribution function of $\m{zq}[V]$ is
obtained as $\Pr\{\m{zq}[V] \le w;v0,tau\} =
\Phi\{\m{zfzq[w,\{q0,q1,q2,q3\},v0,tau]}\}$. We observe that $\m{
zfzq[w,\{0,0,0,0\},0,1]}=w$, and thus the distribution function of z8
under $v0=0$ is $\Pr(Z8\le w;0,1)\approx\Phi(w)$.

\subsubsection{bootstrap probability}

The multiscale bootstrap probability for the scale $\tau$ is \[
\tilde\alpha_1(y,\tau)=\Pr\{V\le0;u0=\hat u,v0=\hat v,\tau\} \] The
corresponding $z$-value is denoted $\tilde
z_1(y,\tau)=-\Phi^{-1}(\tilde\alpha_1(y,\tau)) =-z_c(0;\hat u, \hat v,
\tau)$ in the notation of \refsec{accuracy-bootstrap} of \SH.  For
$y=\eta(\hat u, \hat v)$ with $\hat u=u0=0$ and $\hat v = \m{v}$, and
$\tau=tau1$, the $z$-value is expressed as \[ \m{z1[v,tau1]} =
\m{zc[0,\{0,0,0,0\},v,tau1]} \] in \MN. For $tau1=1$, we define
\m{z0[v]=z1[v,1]}, corresponding to $\hat z_0(y)$ in \SH. For general
value of $tau1$, \m{w1=tau1 z1[v,tau1]} is regarded as another $w$ with
$\m{cb}_r$'s being \m{cbw1} and $\m{c}_r$'s being \m{ccw1}, and the
distribution function is expressed as $\Pr\{W1\le w; v0,tau\} =
\Phi(\m{zfw1})$ under the scale $tau$.  For $tau=1$ and $tau1=1$, we
have $\Pr\{\hat Z_0\le w;v0,tau=1\}=\Phi(\m{zfz0}[w,v0])$, which becomes
$\m{zfz0}[w,0] = w+\oi$ under $v0=0$, showing the first-order accuracy
of \m{z0[v]}.

\subsubsection{double bootstrap}

The z-value of the double bootstrap probability is
\[
 \m{zd}[v] = -\Phi^{-1}( \Pr\{  \hat Z_0 \le \hat z_0(v);v0=0\}),
\]
corresponding to $\hat z_{\rm double}(y)$ in \refsec{accuracy-bootstrap}
of \SH. We observe that $\m{zd[v]}= \m{z8[v]}$, showing the third-order
accuracy of the double bootstrap probability.

\subsubsection{two-level bootstrap}

The ABC formula is given in \m{abcformula[v,ac]}. The $z$-value
corresponding to the the two-level bootstrap corrected $p$-value is
calculated in \m{za[v]}. Its $q_i$'s are in \m{qqza}. The distribution
function of \m{za[v]} under $tau=1$ is \[ \Pr\{\m{za}[V]\le w;v0,1\} =
\Phi(\m{zfzq[w,qqza,v0,1]}), \] which becomes $\Phi(w)+\oii$ for $v0=0$,
showing the second-order accuracy of the two-level bootstrap.

\subsection{Multistep-multiscale bootstrap method}
\subsubsection{a generalization of the pivot}

The pivot \m{z8[v]} is generalized to define
\[
 \m{z8[v,v0,tau]}=\m{zc[v,\{0,0,0,0\},v0,tau]},
\]
which is denoted $\hat z_\infty(0,v;v0,tau)$ in the proof of
\reflemma{threestepboot} of \SH. \m{z8[v,v0,tau]} reduces to \m{z8[v]}
when $v0=0$ and $tau=1$. We show that
\[
 \Pr\{\m{z8}[V,v0,tau]\le w;v0,tau\} \approx \Phi(w).
\]
The inverse function of \m{z8[v,v0,tau]=z} in terms of \m{v} is also
defined here so that \m{v8[z,v0,tau]=v}.

\subsubsection{some useful formula for normal integration}

Let $\m{func}[z]=(az+b)+\m{rem}[z]$, where $\m{rem}[z]$ is a polynomial
function of $z$ with magnitude $\oi$. Then
\[
 \Phi^{-1}\left\{
\int_{-\infty}^\infty \Phi(\m{func}[z])f[z]\,dz
\right\} \approx \m{ddint}[a,b,\m{func2dd}[a,b,\m{rem},z],z].
\]

\subsubsection{two-step multiscale bootstrap}

Using \m{ddint} function defined above, we calculate
\begin{align*}
 \m{z2[v0,tau1,tau2]} &\approx \Phi^{-1}\left\{
\int_{-\infty}^\infty \Phi(\m{z1}[v,tau2]) f[v,v0,tau1]\,dv\right\}\\
 &\approx \Phi^{-1}\left\{
\int_{-\infty}^\infty \Phi(\m{z1}[\m{v8}[z,v0,tau1],tau2])f[z]\,dz
\right\},
 \end{align*}
where \m{z2[v0,tau1,tau2]} corresponds to $\tilde z_2(y,\tau_1,\tau_2)$
in eq.~\refeq{z2-integration} of \SH\ for $y=\eta(\hat u,\hat v)$ with
$\hat u=u0=0$ and $\hat v = v0$. This is the $z$-value of the
twostep-multiscale bootstrap probability.

\subsubsection{three-step multiscale bootstrap}

Similarly, we calculate the $z$-value of the threestep-multiscale
bootstrap probability defined by
\begin{align*}
 \m{z3[v0,tau1,tau2,tau3]} &\approx \Phi^{-1}\left\{
\int_{-\infty}^\infty \Phi(\m{z2}[v,tau2,tau3]) f[v,v0,tau1]\,dv\right\}\\
 &\approx \Phi^{-1}\left\{
\int_{-\infty}^\infty \Phi(\m{z2}[\m{v8}[z,v0,tau1],tau2,tau3])f[z]\,dz
\right\}.
 \end{align*}
This gives the expression of $\tilde z_3(y,\tau_1,\tau_2,\tau_3)$ with
$\hat u=0$ and $\hat v=v0$.

\subsubsection{simplifying z3 and z8}

The six geometric quantities $\gamma_1,\ldots,\gamma_6$ defined in
\reflemma{threestepboot} of \SH\ are denoted by \m{G1,...,G6} here. The
scaling parameters $s_1,\ldots,s_4$ in eq.~\refeq{zeta3} of \SH\ are now
denoted \m{S1,...,S4}. We define \m{Z3G} and \m{Z8G} in terms of
\m{G1,...,G6,S1,...,S4}, and show that
$\m{Z3G}\approx\m{z3[v,tau1,tau2,tau3]}$ and
$\m{Z8G}\approx\m{z8[v]}$. In fact, \m{Z3G} and \m{Z8G}, respectively,
are $\zeta_3(\gamma_1,\ldots,\gamma_6,\tau_1,\tau_2,\tau_3)$ of
eq.~\refeq{zeta3} and $\hat z_3(y)$ of eq.~\refeq{z3} defined in
\refsec{main} of \SH.

We have shown earlier that $\tilde z_3=\m{z3}$ is the $z$-value of the
threestep-multiscale bootstrap probability, and that $\hat
z_\infty=\m{z8}$ is the third-order accurate pivot statistic. Therefore,
the equivalence $\zeta_3 = \m{Z3G} \approx \m{z3} = \tilde z_3$ and
$\hat z_3 = \m{Z8G} \approx \m{z8} = \hat z_\infty$ proves the
third-order accuracy of the threestep-multiscale bootstrap corrected
$p$-value, in which \m{G1,\ldots,G6} are estimated from \m{Z3G} and used
to calculate \m{Z8G}.  Should we note the notational difference between
$\tilde z_3$ and $\hat z_3$; namely, the $z$-value of the
threestep-multiscale bootstrap probability and the $z$-value of the
threestep-multiscale bootstrap corrected $p$-value, respectively.


\end{document}